\documentclass[twocolumn]{autart}    

\usepackage[utf8]{inputenc}
\usepackage[T1]{fontenc}
\usepackage[english]{babel}

\usepackage[retainorgcmds]{IEEEtrantools}
\usepackage{amsfonts,amsmath,amssymb} 

\usepackage{graphicx}
\graphicspath{{Figures/}}

\usepackage{times}
\usepackage{url}
\usepackage{array}
\usepackage{siunitx}

\newcommand{\abs}[1]{\left\lvert#1\right\rvert}

\newcommand{\inner}[2]{\langle#1,#2\rangle}
\newcommand{\SO}{\mathrm{SO(3)}}

\begin{document}

\begin{frontmatter}

\title{A global exponential observer for velocity-aided attitude estimation\thanksref{footnoteinfo}} 

\thanks[footnoteinfo]{This work was supported by the French Agence Nationale de la Recherche through the ANR ASTRID SCAR project ``Sensory Control of Aerial Robots'' (ANR-12-ASTR-0033).
}

\author[Mines]{Philippe Martin}\ead{philippe.martin@mines-paristech.fr},    
\author[Mines]{Ioannis Sarras}\ead{sarras@ieee.org},              
\author[I3S]{Minh-Duc Hua}\ead{hua@i3s.unice.fr},  
\author[I3S]{Tarek Hamel}\ead{thamel@i3s.unice.fr}

\address[Mines]{Centre Automatique et Systèmes, Mines ParisTech, PSL Research University, Paris, France}  
\address[I3S]{I3S UNS-CNRS, Sophia-Antipolis, France}        


\begin{abstract}                          
We propose a simple nonlinear observer for estimating the attitude and velocity of a rigid body from the measurements of specific acceleration, angular velocity, magnetic field (in body axes), and linear velocity (in body axes). It is uniformly globally exponentially convergent, and also enjoys other nice properties: global decoupling of pitch and roll estimation from magnetic measurements, good local behavior, and easy tuning. The observer is ``geometry-free'', in the sense that it respects only asymptotically the rotational geometry. The good behavior of the observer, even when the measurements are noisy and biased is illustrated in simulation.
\end{abstract}

\end{frontmatter}

\section{Introduction}\label{sec:intro}
Estimating the attitude and velocity of a robotic vehicle is usually crucial for control purposes. For lightweight and low-cost systems equipped with a MEMS Inertial Measurement Unit, ``true'' inertial navigation (i.e., based on the Schuler effect due to a rotating non-flat Earth) is excluded, because such sensors are not accurate enough in the long run. To estimate in particular the velocity, the inertial sensors must be ``aided''. This can be done thanks to a model of the forces acting on the vehicle, see e.g.~\cite{MartiS2010ICRA,MartiS2015arXiv} in the context of quadrotors, in which case the estimator is specific to the vehicle. The alternative is to use an additional velocity sensor, providing the velocity vector in Earth axes, see e.g.~\cite{MartiS2008IFAC,Hua2010CEP,RoberT2011CDC,GripFJS2012ACC}, or in body axes, see e.g.~\cite{BonnaMR2008ITAC,DukanS2013CCAMS,TroniWICRA2013,HuaMH2016AUT}; in this case the estimator can be generic, but at the cost of the extra sensor. In particular, two nonlinear observers are proposed in~\cite{HuaMH2016AUT} for the estimation with velocity aiding in body axes; these observers enjoy nice properties: almost global asymptotic stability, global decoupling of roll and pitch estimation from magnetic measurements, good local behavior and easy tuning. Notice by the way that it is quite difficult to achieve the same level of performance with an estimator based on an Extended Kalman Filter (or its variants).

The present paper is in some sense a sequel to~\cite{HuaMH2016AUT}: it considers exactly the same problem, namely attitude and velocity estimation from rate gyros, accelerometers, magnetometers, and body-axes velocity. The proposed observer has the same desirable features as those in~\cite{HuaMH2016AUT}, but is \emph{uniformly globally exponentially stable} instead of merely \emph{almost globally asymptotically stable}; moreover, it is very simple, and so is the convergence proof. To achieve this result, the idea is to ``forget'' the geometry of~$\SO$ (or of the unit quaternion space) on which the orientation lives: whereas the observers of~\cite{HuaMH2016AUT}, which are instances of so-called invariant observers~\cite{BonnaMR2008ITAC}, respect the rotational symmetries at each time, the present observer lives in a bigger space and respects them only asymptotically. This idea of designing an observer on a bigger space relaxing the geometric constraints of the model is quite recent, see e.g.~\cite{BatisSO2012CEP,BatisSO2012AUT,EudesM2014IROS,BatisSO2014CDC,GripFJS2015AUT,MartiS2015arXiv,MartiS2016arXiv}. 

The paper runs as follows: section~\ref{sec:model} introduces the design model (dynamics and measurements) on which the observer is  based; section~\ref{sec:observer} presents the observer, and studies its convergence; finally, section~\ref{sec:simulation} illustrates in simulation the good behavior of the observer, even when the measurements are noisy and biased.

\section{The design model}\label{sec:model}
We consider a moving rigid body subjected to the angular velocity~$\omega$ (in body axes). Its orientation (from inertial to body axes) matrix $R\in\SO$ is related to~$\omega$ by
\begin{IEEEeqnarray}{rCl}
	\dot R &=& RS(\omega)\label{eq:R},
\end{IEEEeqnarray} 
where $S(\omega)$ is the skew-symmetric matrix defined by $S(\omega)x:=\omega\times x$ whatever the vector~$x$.
Assuming a flat and non-rotating Earth, the specific acceleration $a$ (in body axes) of a point of the body is by definition related to the velocity of this point ($V$ in Earth axes, $v$ in body axes) by
\begin{IEEEeqnarray}{rCl}
	a &=& R^T(\dot V-gE_3)=\dot v+\omega\times v-gR^TE_3\label{eq:a},
\end{IEEEeqnarray} 
where $gE_3$ is the gravity vector (in Earth axes).
Setting $\gamma:=gR^TE_3$, we obviously have $\dot{\gamma}=\gamma\times\omega$ since $E_3$ is constant; similarly, if $B$ is some other constant vector (in Earth axes), $\beta:=R^TB$ satisfies $\dot{\beta}=\beta\times\omega$. Provided $B$ and $E_3$ are not collinear, these two differential equations for $\gamma$ and $\beta$ are equivalent to~\eqref{eq:R}; indeed the rotation matrix~$R$ is completely specified by $\gamma$ and~$\beta$, since
\begin{IEEEeqnarray}{rCl}
	R^T &=& R^T\begin{pmatrix}gE_3& B& gE_3\times B\end{pmatrix}\cdot\begin{pmatrix}gE_3& B& gE_3\times B\end{pmatrix}^{-1}\nonumber\\
	&=& \begin{pmatrix}\gamma& \beta& \gamma\times\beta\end{pmatrix}\cdot\begin{pmatrix}gE_3& B& gE_3\times B\end{pmatrix}^{-1},
\end{IEEEeqnarray} 
where we have used $R^T(x\times y)=R^Tx\times R^Ty$ whatever the vectors $x,y$, since $R$ is a rotation matrix.
The dynamics~\eqref{eq:R}-\eqref{eq:a} of the rigid body can therefore be expressed as
\begin{IEEEeqnarray}{rCl}
	\dot{v} &=& v\times\omega + \gamma + a \label{eq:v}\\
	\dot{\gamma} &=& \gamma\times\omega \label{eq:gamma}\\
	\dot{\beta} &=& \beta\times\omega \label{eq:beta}.
\end{IEEEeqnarray} 
Notice the two subsystems \eqref{eq:v}-\eqref{eq:gamma} and~\eqref{eq:beta} are completely independent, a property which is hidden in~\eqref{eq:R}-\eqref{eq:a}.

We assume that the rigid body is equipped with a 3-axis linear velocity sensor, e.g.~a Doppler radar, together with a strapdown unit comprising a 3-axis rate gyro, a 3-axis accelerometer, and a 3-axis magnetometer. If these sensors were ideal, they would provide the perfect measurements
\begin{IEEEeqnarray}{rCl}
	v_m &=& v\label{eq:mesv}\\
	\omega_m &=& \omega\\
	a_m &=& a\\
	\beta_m &=& \beta\label{eq:mesbeta}.
\end{IEEEeqnarray} 
Notice that the magnetic field~$B$ (in Earth axes) is (locally) constant, and not collinear with $E_3$ (except at the Earth magnetic poles).
Of course, the sensors are in fact corrupted by biases and noises, and actually provide the measurements
\begin{IEEEeqnarray}{rCl}
	v_m &=& v+b_v+\nu_v\label{eq:mesva}\\
	\omega_m &=& \omega+b_\omega+\nu_\omega\\
	a_m &=& a+b_a+\nu_a\\
	\beta_m &=& \beta+b_\beta+\nu_\beta,\label{eq:mesbetaa}
\end{IEEEeqnarray} 
where the $b_i$, $i=v,\omega,a,\beta$, are constant (or slowly-varying) biases, and the $\nu_i$ are (more or less Gaussian white) noises.

The design model on which the observer is based and its convergence analyzed consists of the dynamics~\eqref{eq:v}--\eqref{eq:beta} with perfect measurements~\eqref{eq:mesv}--\eqref{eq:mesbeta}; the behavior when using the actual measurements~\eqref{eq:mesva}--\eqref{eq:mesbetaa} will be tested in simulation.
\section{The observer and its analysis}\label{sec:observer}
\subsection{The observer}
We are going to show that the state of~\eqref{eq:v}--\eqref{eq:beta} can be estimated by the observer
\begin{IEEEeqnarray}{rCl}
	\dot{\hat{v}} &=& \hat{v}\times\omega_m + a_m + \hat{\gamma} - (L+K)(\hat{v}-v_m)\label{eq:hatv}\\
	\dot{\hat{\gamma}} &=& \hat{\gamma}\times\omega_m - \bigl(LS(\omega_m) - S(\omega_m)L + LK\bigr)(\hat{v}-v_m)\IEEEeqnarraynumspace\label{eq:hatgamma}\\
	\dot{\hat{\beta}} &=& \hat{\beta}\times\omega_m - M(\hat{\beta}-\beta_m),\label{eq:hatbeta}
\end{IEEEeqnarray}
where the $3\times3$ matrices $K,L,M$ are tuning parameters. Notice this observer is a copy of~\eqref{eq:v}--\eqref{eq:beta} with (time-varying) correction terms, which respects the decoupled structure of~\eqref{eq:v}--\eqref{eq:beta}.
The independence of $\hat\gamma$ from the easily perturbed magnetic measurements is a very desirable feature. Indeed, if the orientation matrix~$R$ is parametrized by the roll, pitch and yaw Euler angles~$(\phi,\theta,\psi)$, $\gamma$ reads
\begin{IEEEeqnarray*}{rCl}
	\gamma &=& g\begin{pmatrix}-\sin\theta& \sin\phi\cos\theta& \cos\phi\cos\theta\end{pmatrix}^T;
\end{IEEEeqnarray*} 
in other words, $\gamma$ encodes the roll and pitch angles, which are crucial for stabilization purposes (whereas the yaw angle far less matters).

Also notice the observer was derived using the invariant-manifold-observer methodology proposed in~\cite{InIbook,KaragSA2009AUT}; this approach, seemingly well-adapted to problems in aerial robotics, has already been used in this context in e.g.~\cite{SarraS2014MSC,MartiS2016arXiv}.

\subsection{Convergence analysis}
Defining the error variables
\begin{IEEEeqnarray*}{rCl}
	e_v &:=& \hat{v} - v\\
	e_\gamma &:=& \hat{\gamma} - \gamma - Le_v\\
	e_\beta &:=& \hat{\beta} - \beta,
\end{IEEEeqnarray*} 
and assuming perfect measurements, the error system reads
\begin{IEEEeqnarray}{rCl}
	\dot{e}_v &=& e_v\times\omega + e_\gamma - K e_v \label{eq:ev}\\
	\dot{e}_\gamma &=& e_\gamma\times\omega - Le_\gamma \\
	\dot{e}_\beta &=&  e_\beta\times\omega - Me_\beta \label{eq:ebeta}.
\end{IEEEeqnarray} 

\begin{thm}\label{th:UGES}
	If the symmetric parts of $K,L,M$ are positive definite, the equilibrium point $(\bar e_v,\bar e_\gamma,\bar e_\beta):=(0,0,0)$ of the error system~\eqref{eq:ev}--\eqref{eq:ebeta} is globally exponentially stable.
\end{thm}

\begin{pf}
Consider the candidate Lyapunov function
\begin{IEEEeqnarray*}{rCl}
	V(e_\gamma,e_v,e_\beta)&:=& \underbrace{\frac{\rho_1}{2}\abs{e_\gamma}^2 + \frac{\rho_1\epsilon^2}{2}\abs{e_v}^2}_{=:V_\gamma(e_\gamma,e_v)}
	+\underbrace{\frac{1}{2}\abs{e_\beta}^2}_{=:V_\beta(e_\beta)},
\end{IEEEeqnarray*} 
with $\rho_1,\epsilon>0$. On the one hand,
\begin{IEEEeqnarray*}{rCl}
	\dot V_\beta(e_\beta) &=& -\inner{e_\beta}{\frac{M+M^T}{2}e_\beta} \le -\underline{\sigma}_M\abs{e_\beta}^2,
\end{IEEEeqnarray*} 
where we have used $\inner{x}{x\times y}=0$, and denoted by $\underline{\sigma}_M$ the smallest eigenvalue of~$\frac{M+M^T}{2}$, which is strictly positive by assumption. On the other hand,
\begin{IEEEeqnarray*}{rCl}
	\dot V_\gamma &=& -\rho_1\inner{e_\gamma}{\frac{L+L^T}{2}e_\gamma}+\rho_1\epsilon^2\inner{e_v}{e_\gamma}\\
	&&-\,\rho_1\epsilon^2\inner{e_v}{\frac{K+K^T}{2} e_v}\\
	&\leq& -\rho_1\underline{\sigma}_L\abs{e_\gamma}^2 
	+ \rho_1\epsilon^2\biggl(\frac{\epsilon\abs{e_v}^2}{2} + \frac{\abs{e_\gamma}^2}{2\epsilon}\biggr) -\rho_2\underline{\sigma}_K\abs{e_v}^2\\
	&=& -\rho_1\Bigl(\underline{\sigma}_L-\frac{\epsilon}{2}\Bigr)\abs{e_\gamma}^2 
	- \rho_1\epsilon^2\Bigl(\underline{\sigma}_K-\frac{\epsilon}{2}\Bigr)\abs{e_v}^2,
\end{IEEEeqnarray*} 
where Young's inequality $\inner{x}{y}\leq \frac{\epsilon\abs{x}^2}{2}+\frac{\abs{y}^2}{2\epsilon}$ has been applied to the cross term.
Since $\epsilon$ can be chosen as small as desired and $\underline{\sigma}_L,\underline{\sigma}_K>0$ by assumption, $V$ is clearly a strict Lyapunov function, which proves the claim.\qed
\end{pf}

There remains to build an estimate of the orientation matrix~$R$ from the estimated vectors $\hat\gamma$ and~$\hat\beta$. It is convenient to choose for Earth axes
the North-East-Down frame $(E_1,E_2,E_3)$, in which the magnetic vector reads $B=(B_1,0,B_3)^T$. We then have the obvious but important following corollary.

\begin{cor}
	Under the assumptions of Theorem~\ref{th:UGES},
	\begin{IEEEeqnarray*}{rCl}
		\tilde R &:=& \begin{pmatrix}\dfrac{(\hat\gamma\times\hat\beta)\times\hat\gamma}{\abs{(gE_3\times B)\times gE_3}}& 
		\dfrac{\hat\gamma\times\hat\beta}{\abs{gE_3\times B}} & 
		\dfrac{\hat\gamma}{\abs{gE_3}}\end{pmatrix}^T
	\end{IEEEeqnarray*} 
	globally exponentially converges to the orientation matrix~$R$.
\end{cor}
\begin{pf}
	By Theorem~\ref{th:UGES}, $\hat\gamma\to\gamma$ and $\hat{\beta}\to\beta$. Hence,
	\begin{IEEEeqnarray*}{rCl}
		\tilde R &\to& \begin{pmatrix}\frac{(\gamma\times\beta)\times\gamma}{\abs{(gE_3\times B)\times gE_3}}& 
			\frac{\gamma\times\beta}{\abs{gE_3\times B}} & 
			\frac{\gamma}{\abs{gE_3}}\end{pmatrix}^T\\
		&=& \begin{pmatrix}\frac{(gE_3\times B)\times gE_3}{\abs{(gE_3\times B)\times gE_3}}& 
			\frac{gE_3\times B}{\abs{gE_3\times B}} & 
			\frac{gE_3}{\abs{gE_3}}\end{pmatrix}^T\cdot R\\
		&=& R,
	\end{IEEEeqnarray*} 
where we have used $R^T(x\times y)=R^Tx\times R^Ty$, and $gE_3\times B=gB_1E_2$ and $(gE_3\times B)\times gE_3=g^2B_1E_1$.\qed
\end{pf}

Of course, $\tilde R$ has no reason to be a rotation matrix (it is only asymptotically so); it has nevertheless orthogonal (possibly zero) rows. If a bona fide rotation matrix is required at all times, a natural idea is to project $\tilde R$ on the ``closest'' rotation matrix~$\hat R$, thanks to a polar decomposition. Because $\tilde R$ has orthogonal rows, the expression of $\hat R$ is readily found, without using the standard but computationally heavy projection algorithm based on singular value decomposition. For details about the polar decomposition and related matters, see e.g.~\cite[Chapter~$8$]{Higham2008book}.

\begin{prop}
	Considerer the polar decomposition of $\tilde R^T$
	\begin{IEEEeqnarray*}{rCl}
		\tilde R^T &=& \hat R^T(\tilde R\tilde R^T)^\frac{1}{2}.
	\end{IEEEeqnarray*} 
	Then $\hat R$, which is by construction the best approximation of~$\tilde R$ among all orthogonal matrices, is a rotation matrix that globally exponentially converges to~$R$. When $\hat\gamma$ and $\hat\beta$ are not collinear, $\hat R$ is uniquely defined by
	\begin{IEEEeqnarray*}{rCl}
		\hat R &:=& \begin{pmatrix}\dfrac{(\hat\gamma\times\hat\beta)\times\hat\gamma}{\abs{(\hat\gamma\times\hat\beta)\times\hat\gamma}}& 
			\dfrac{\hat\gamma\times\hat\beta}{\abs{\hat\gamma\times\hat\beta}} & 
			\dfrac{\hat\gamma}{\abs{\hat\gamma}}\end{pmatrix}^T.
	\end{IEEEeqnarray*} 
\end{prop}
\begin{pf}
	Since $\tilde R$ has orthogonal rows,
	\begin{IEEEeqnarray*}{rCl}
		(\tilde R\tilde R^T)^\frac{1}{2}
		&=& \begin{pmatrix}\frac{\abs{(\hat\gamma\times\hat\beta)\times\hat\gamma}}{\abs{(gE_3\times B)\times gE_3}}& 0&0\\
			0& \frac{\abs{\hat\gamma\times\hat\beta}}{\abs{gE_3\times B}} & 0\\
		0&0& \frac{\abs{\hat\gamma}}{\abs{gE_3}}\end{pmatrix}.
	\end{IEEEeqnarray*} 
	When $\hat\gamma$ and $\hat\beta$ are not collinear, the expression for $\hat R$ follows at once from $\hat R^T=\tilde R^T(\tilde R\tilde R^T)^{-\frac{1}{2}}$. When $\hat\gamma=0$, one may choose $\hat R:=I$; when $\hat\gamma\neq0$ but $\hat\gamma\times\hat\beta=0$, one may choose $\hat R:=(\hat E_1,\hat E_2,\frac{\hat\gamma}{\abs{\hat\gamma}})^T$, where $\hat E_1$, $\hat E_2$ and $\frac{\hat\gamma}{\abs{\hat\gamma}}$ form a direct orthonormal frame.\qed
\end{pf}

Notice the knowledge of the magnetic vector~$B$ is not used in the observer itself; it is only required for the reconstruction of the full estimated orientation $\tilde R$ or~$\hat R$. It is not even necessary for reconstructing the roll and pitch angles.

\subsection{Gain tuning and local behavior}
Global convergence is certainly a desirable property for an observer: it ensures a ``reasonable'' behavior under exceptional circumstances when the estimated state is far from the actual state of the system; however, the local behavior around ``interesting'' trajectories (e.g. at least nominal equilibrium points) is in practice also of paramount importance.

The observers of~\cite{HuaMH2016AUT} have in that regard an interesting feature, inherited from their invariance properties~\cite{BonnaMR2008ITAC}: their error systems expressed in suitable coordinates are autonomous, i.e., do not depend on the trajectory followed by the system,  whatever the tuning gains. Thanks to this property, the local behavior of such observers is very easy to understand, which renders the tuning simple; indeed, the local behavior of the error system  around \emph{every} trajectory of the body is entirely ruled by the eigenvalues of its tangent linearization (which is time-invariant). The proposed observer does not in general enjoy this property, but is nonetheless very easy to tune. It is convenient to express the error system~\eqref{eq:ev}--\eqref{eq:ebeta} in the rotated coordinates  $(E_v,E_\gamma,E_\beta):=(Re_v,Re_\gamma,Re_\beta)$, which are reminiscent of the invariant error coordinates used in~\cite{HuaMH2016AUT}. This yields
\begin{IEEEeqnarray}{rCl}
	\dot{E}_v &=& E_\gamma - RKR^TE_v \label{eq:Ev}\\
	\dot{E}_\gamma &=& -RLR^TE_\gamma \label{eq:Egamma}\\
	\dot{E}_\beta &=&  - RMR^TE_\beta, \label{eq:Ebeta}
\end{IEEEeqnarray} 
which is a linear system, albeit a priori time-varying because of the presence of the orientation~$R$ of the body. 
The eigenvalues of the subsystem~\eqref{eq:Ebeta} are simply those of the matrix~$RMR^T$; the eigenvalues of the subsystem~\eqref{eq:Ev}-\eqref{eq:Egamma}, thanks to the cascade structure, are those of the matrices $RKR^T$ and~$RLR^T$; of course, the eigenvalues do not characterize stability for a time-varying system.

An obvious choice of the tuning matrices satisfying the assumptions of theorem~\ref{th:UGES} is $(K,L,M):=(kI,lI,mI)$, with $k,l,m$ strictly positive numbers. In this case, \eqref{eq:Ev}--\eqref{eq:Ebeta} is time-invariant and reads
\begin{IEEEeqnarray}{rCl}
	\dot{E}_v &=& E_\gamma - kE_v\label{eq:EvId}\\
	\dot{E}_\gamma &=& -lE_\gamma\label{eq:EgammaId}\\
	\dot{E}_\beta &=&  - mE_\beta; \label{eq:EbetaId}
\end{IEEEeqnarray} 
moreover, \eqref{eq:EvId}-\eqref{eq:EgammaId} splits component-wise into 3 identical subsystems with eigenvalues~$-k,-l$, and \eqref{eq:EbetaId} into 3 identical subsystems with eigenvalue~$-m$. The behavior of the error system, hence the tuning, is therefore very simple. Notice also the simple form of the observer itself, which becomes
\begin{IEEEeqnarray*}{rCl}
	\dot{\hat{v}} &=& \hat{v}\times\omega_m + a_m + \hat{\gamma} - (l+k)(\hat{v}-v_m)\\
	\dot{\hat{\gamma}} &=& \hat{\gamma}\times\omega_m - lk(\hat{v}-v_m)\\
	\dot{\hat{\beta}} &=& \hat{\beta}\times\omega_m - m(\hat{\beta}-\beta_m).
\end{IEEEeqnarray*}

In some cases, it may not be desirable to have as in the previous tuning the same gains on the components of each correction term, for instance when the velocity sensor has for technological reasons a ``privileged'' direction (in general the vertical axis). A possible tuning for \eqref{eq:Ev}--\eqref{eq:Ebeta} is then
\begin{IEEEeqnarray*}{rCl'rCl}
	K &:=& \begin{pmatrix}k_x& -k_y& 0\\ k_y& k_x& 0\\ 0& 0& k_z\end{pmatrix},
	& L &:=& \begin{pmatrix}l_x& -l_y& 0\\ l_y& l_x& 0\\ 0& 0& l_z\end{pmatrix},
\end{IEEEeqnarray*} 
with $k_x,k_z,l_x,l_z>0$; notice~$K$ (resp.~$L$) has a pair of complex conjugate eigenvalues when $k_y\neq0$ (resp.~$l_y\neq0$). For a trajectory of the system such that
\begin{IEEEeqnarray*}{rCl}
	R &:=& \begin{pmatrix}\cos\psi & -\sin\psi& 0\\ \sin\psi& \cos\psi& 0\\ 0& 0& 1\end{pmatrix},
\end{IEEEeqnarray*} 
where $\psi$ is an arbitrary function of time, it is easy to check that $RKR^T=K$ and $RLR^T=L$. On such a trajectory (and approximately so on nearby trajectories), the error subsystem~\eqref{eq:Ev}-\eqref{eq:Egamma} reads
\begin{IEEEeqnarray*}{rCl}
	\dot{E}_v &=& E_\gamma - KE_v\\
	\dot{E}_\gamma &=& - LE_\gamma,
\end{IEEEeqnarray*} 
hence is time-invariant. Notice this situation, which corresponds to the body moving level with an arbitrary velocity while spinning around a vertical axis, approximately corresponds to ``normal'' operation of a vehicle when it is not aggressively maneuvering.

If one is interested only in trajectories where $R\approx I$ (i.e., the body is moving level without spinning), the tuning is very flexible: it is possible to have for $K,L,M$, hence for~\eqref{eq:Ev}--\eqref{eq:Ebeta} with $R:=I$, any eigenvalues with negative real parts (of course complex eigenvalues must come by conjugate pairs); the only excluded configuration is three pairs of complex conjugate eigenvalues for the subsystem~\eqref{eq:Ev}--\eqref{eq:Egamma}.

\section{Simulations}\label{sec:simulation}
\begin{figure}[ht]\hspace{-5mm}
	\includegraphics[width=1.1\columnwidth]{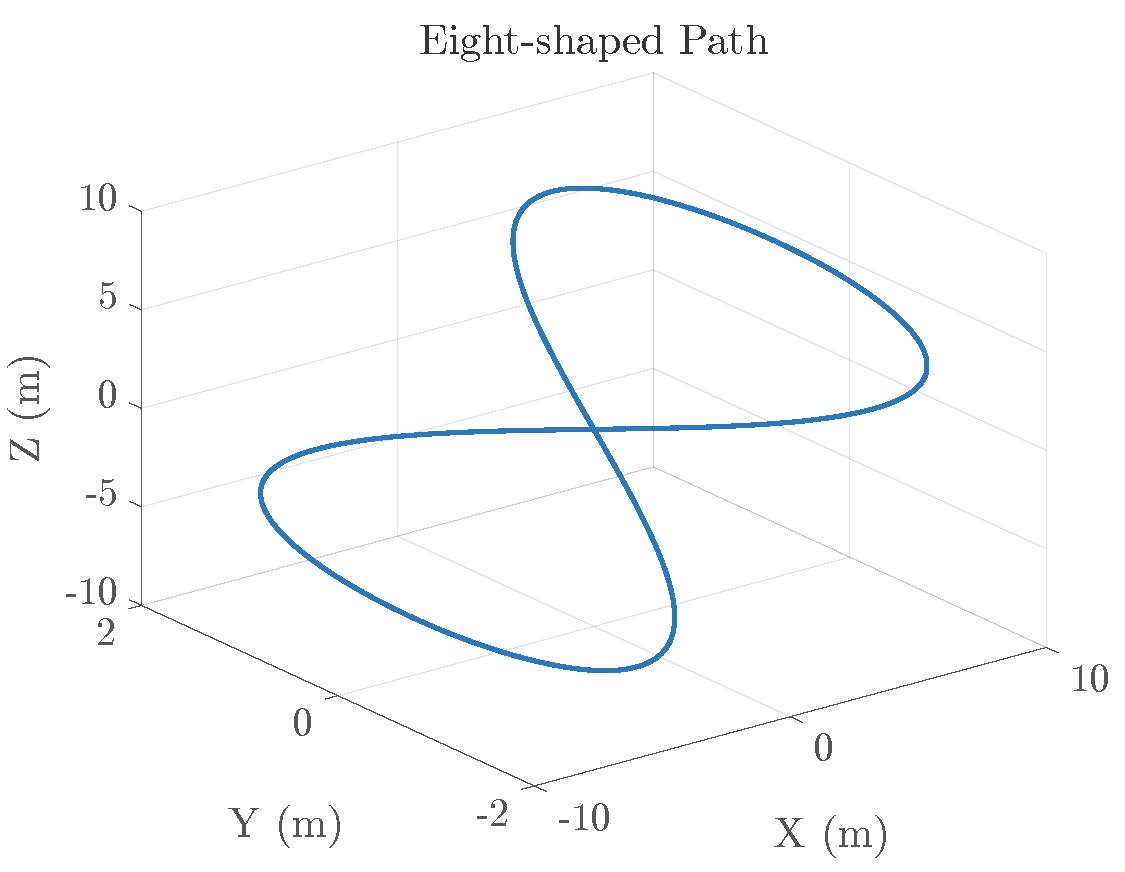}
	\caption{Path followed by the body.}
	\label{fig:path}
\end{figure}

\begin{figure}[ht]\hspace{-5mm}
	\includegraphics[width=1.1\columnwidth]{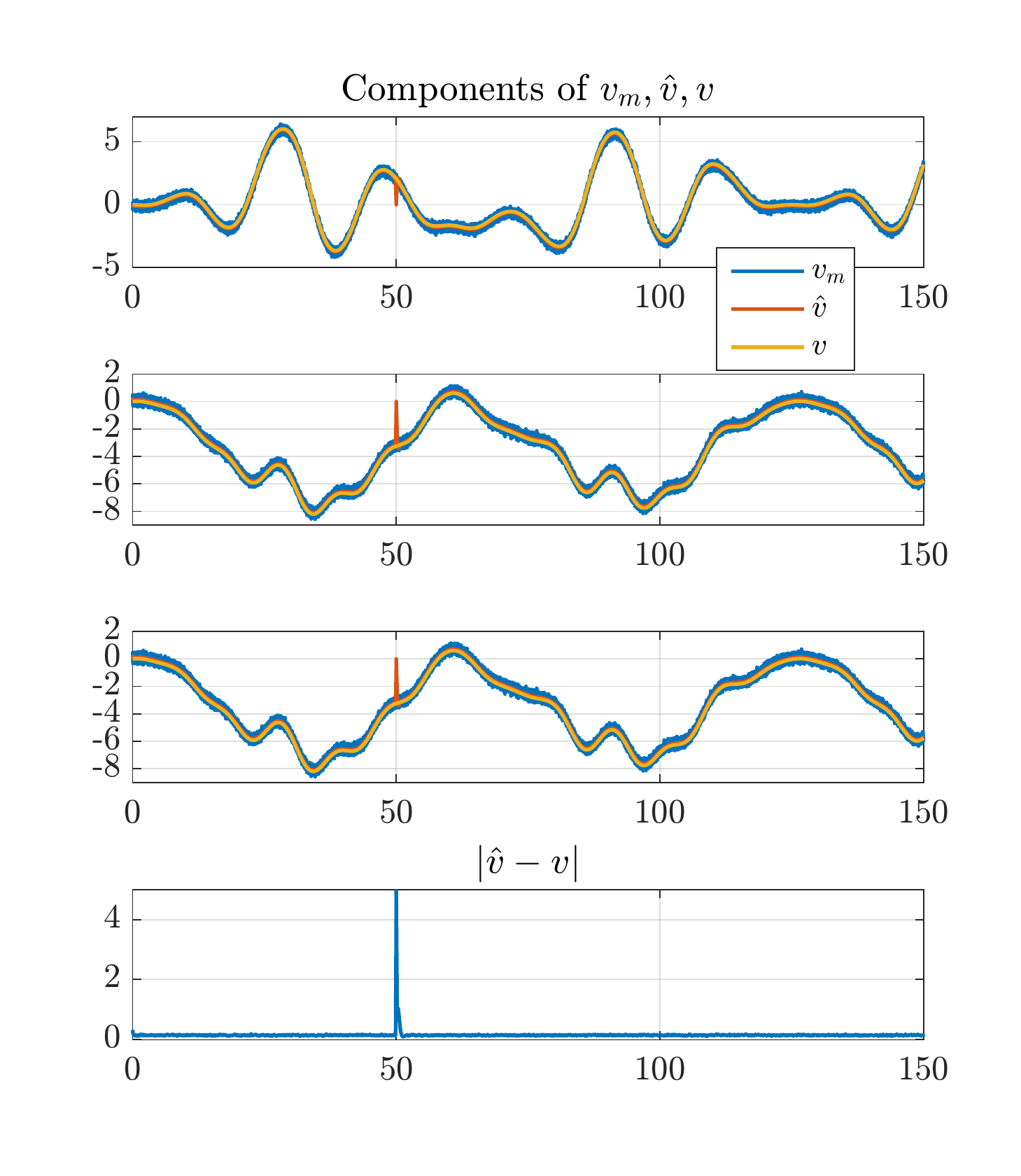}
	\caption{Components of true $v$ (red), measured $v_m$ (blue) and estimated~$\hat v$ (orange).}
	\label{fig:v}
\end{figure}

\begin{figure}[ht]\hspace{-5mm}
	\includegraphics[width=1.1\columnwidth]{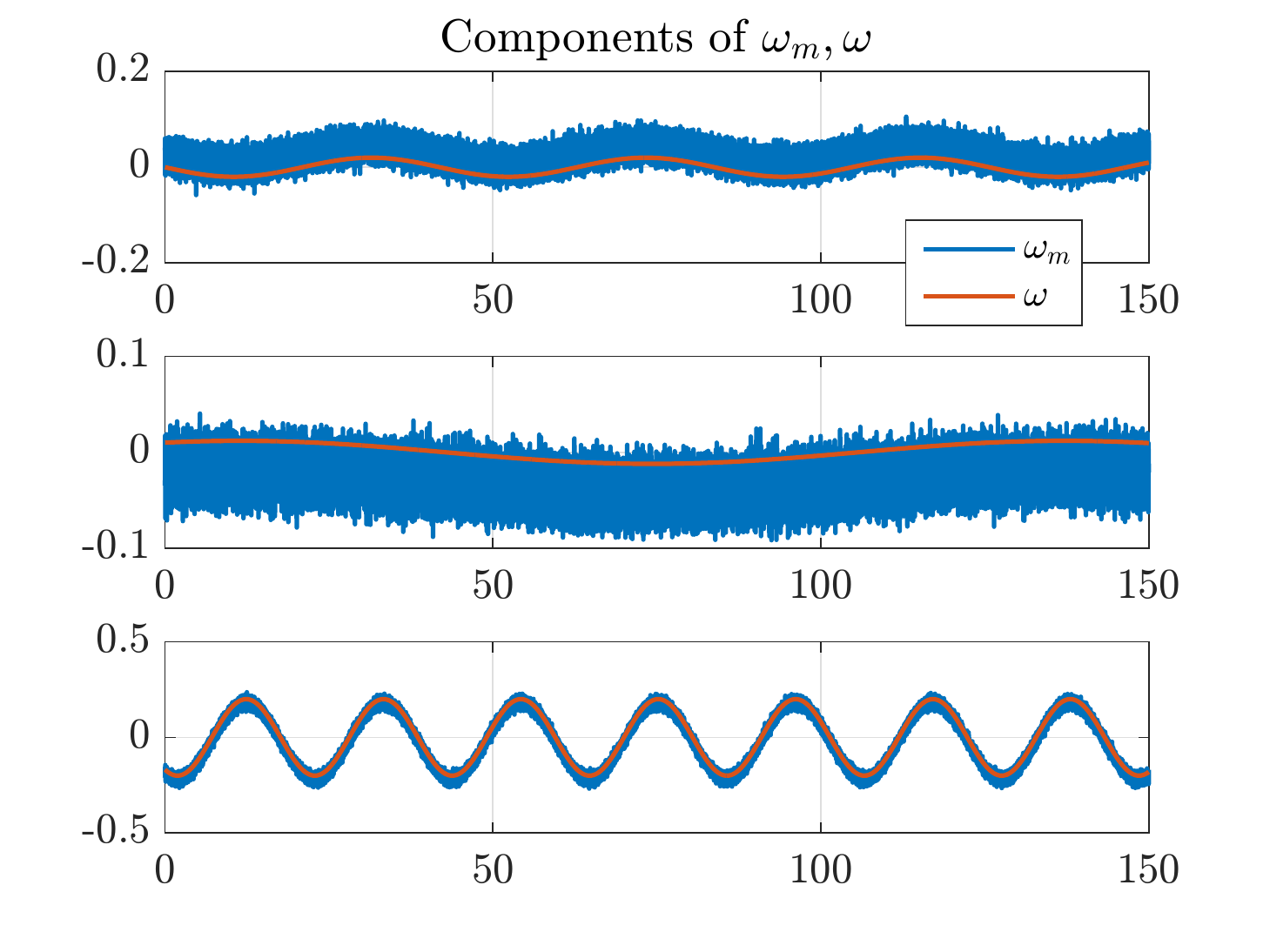}
	\caption{Components of true $\omega$ (red) and measured $\omega_m$ (blue).}
	\label{fig:omega}
\end{figure}

\begin{figure}[ht]\hspace{-5mm}
	\includegraphics[width=1.1\columnwidth]{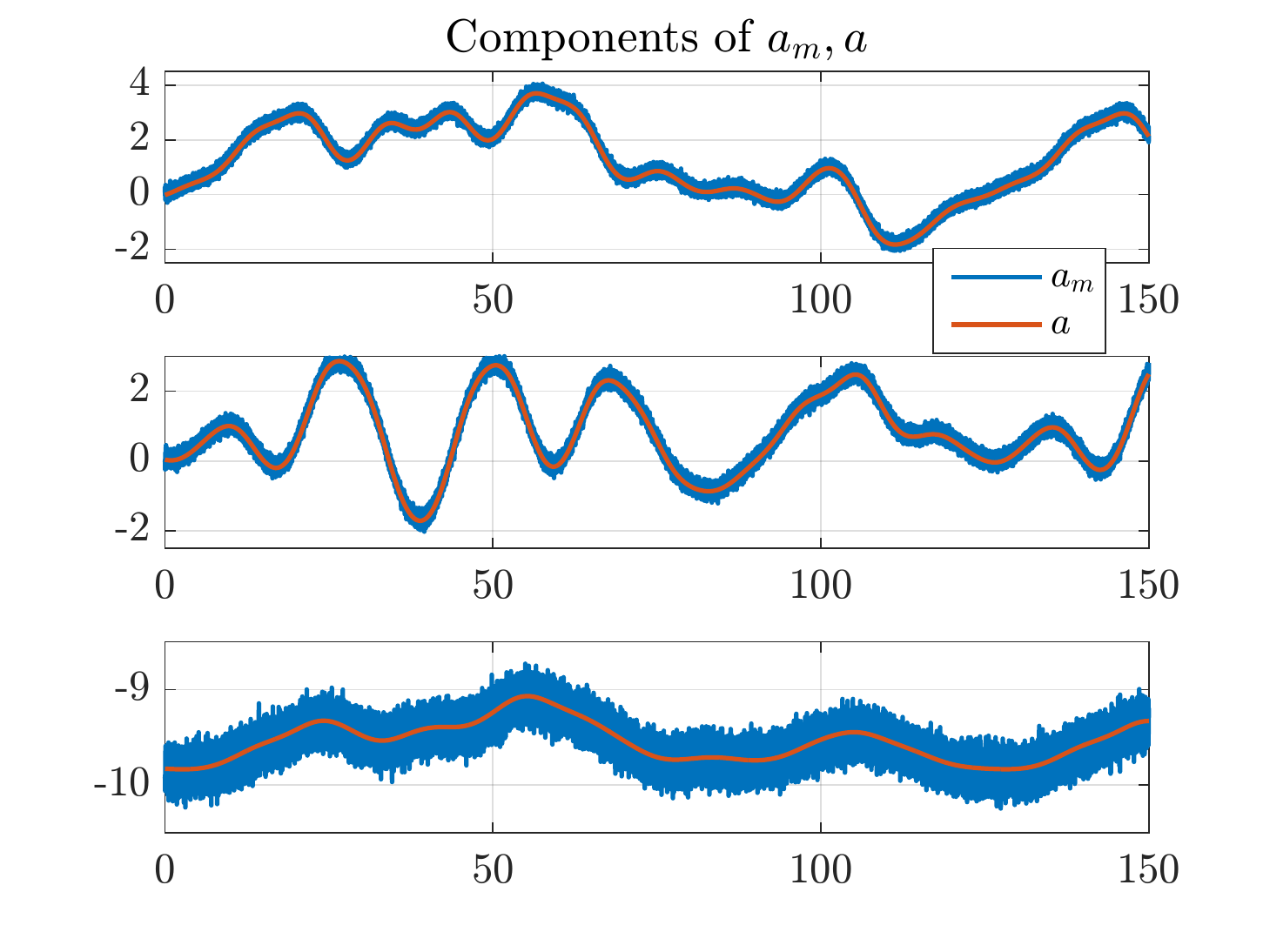}
	\caption{Components of true $a$ (red) and measured $a_m$ (blue).}
	\label{fig:a}
\end{figure}

\begin{figure}[ht]\hspace{-5mm}
	\includegraphics[width=1.1\columnwidth]{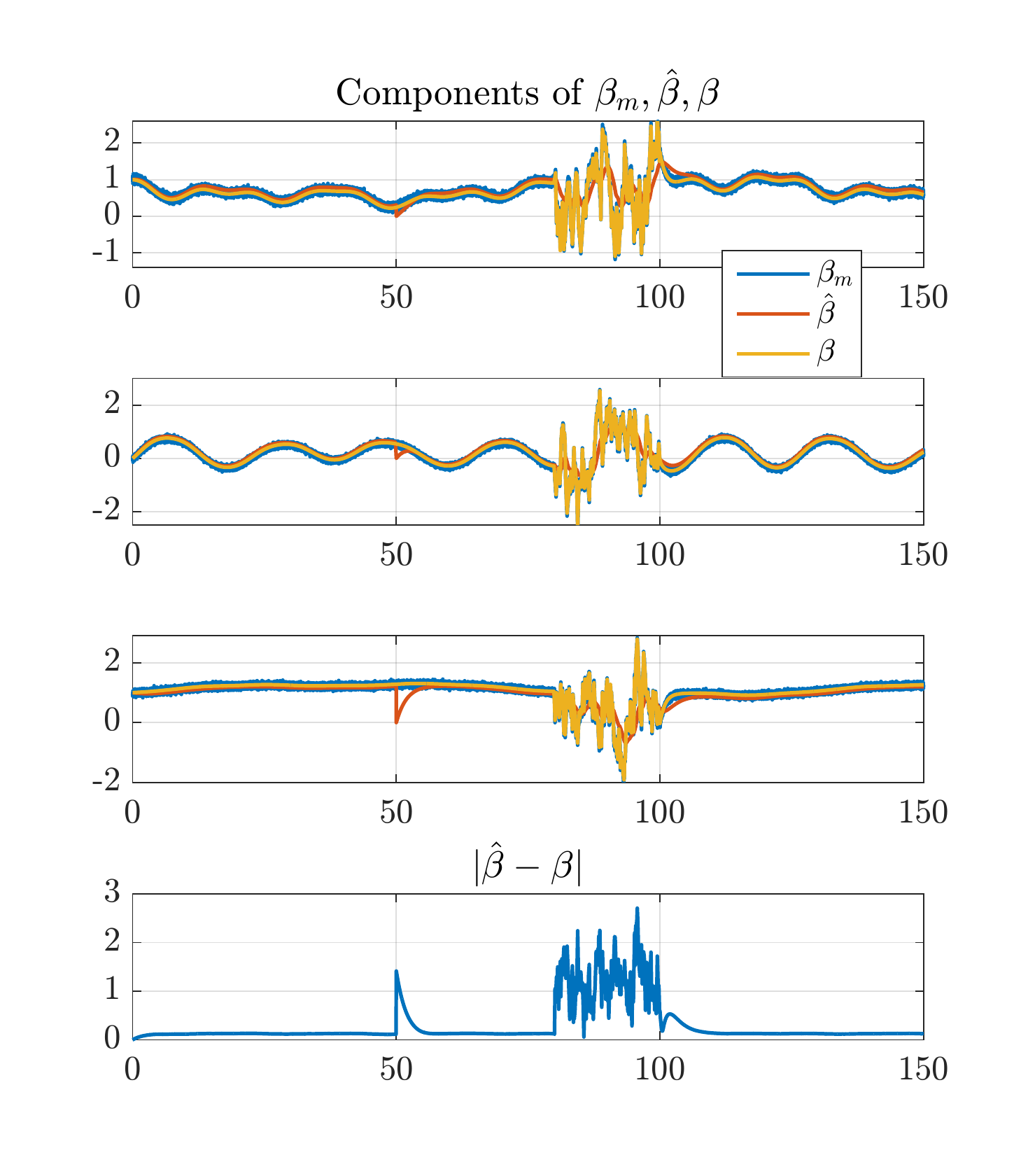}
	\caption{Components of true~$\beta$ (red), measured~$\beta_m$ (blue) and estimated~$\hat\beta$ (orange).}
	\label{fig:beta}
\end{figure}

\begin{figure}[ht]\hspace{-5mm}
	\includegraphics[width=1.1\columnwidth]{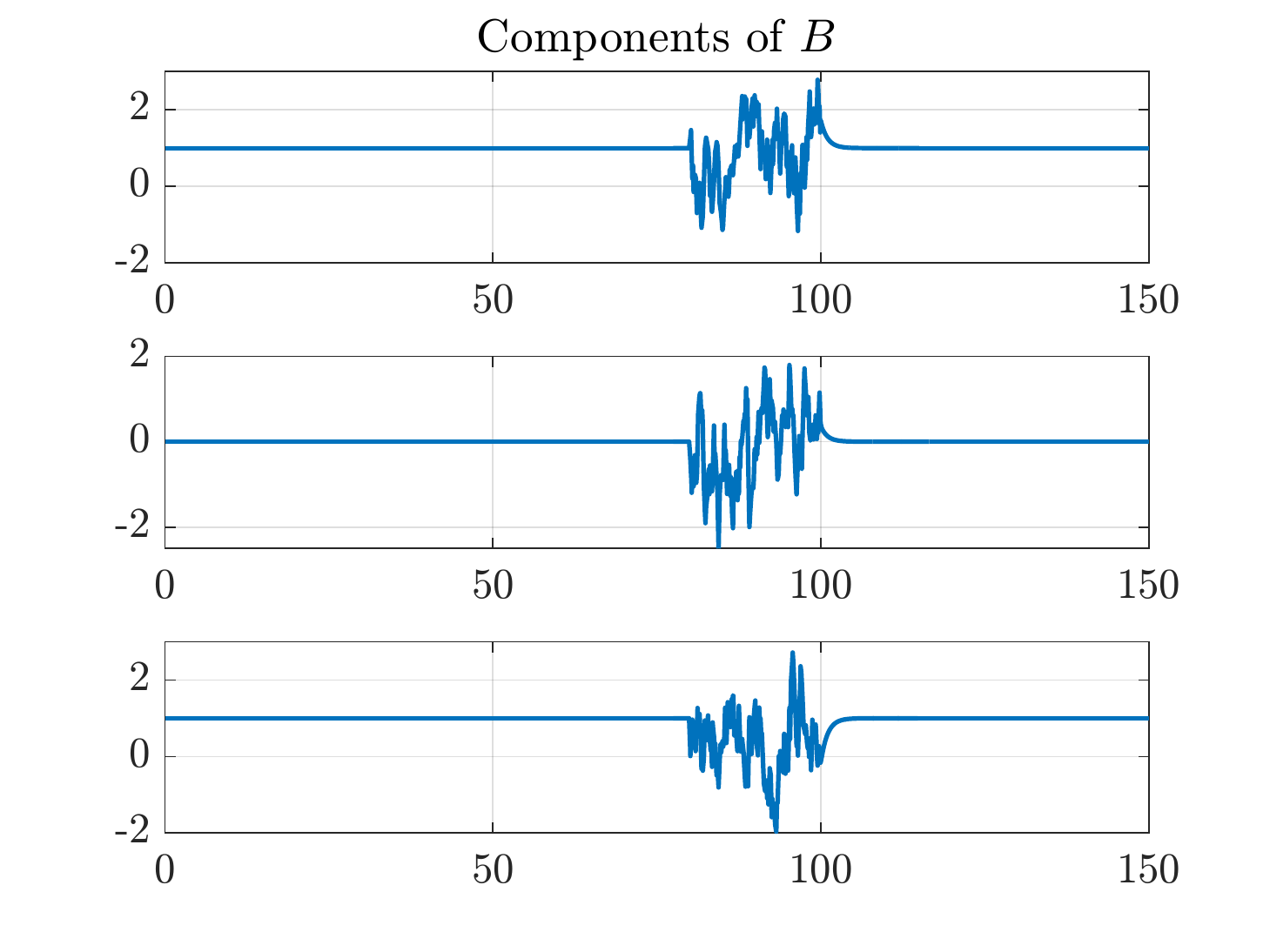}
	\caption{Components of the ``constant'' magnetic vector~$B$.}
	\label{fig:B}
\end{figure}

The good behavior of the observer is now illustrated in a simulation including effects not taken into account in the design model.

The scenario is the following: the body 
follows a tilted ``eight-shaped'' path at the linear velocity~$v(t)$ (and corresponding acceleration~$a(t)$), while undergoing the angular velocity~$\omega(t)$, see Fig.~\ref{fig:path}--\ref{fig:omega}. The magnetic vector~$B$ is set to the nominal value $(1/\sqrt2,0,1/\sqrt2)^T$, but is subjected to a violent disturbance for $t\in[80,100]$, see Fig.~\ref{fig:B}; notice~\eqref{eq:beta} holds only approximately during the disturbance. Finally, $g:=9.81$.

The observer is fed with the imperfect measurements~\eqref{eq:mesva}--\eqref{eq:mesbetaa}, see~Fig.~\ref{eq:v}-\ref{fig:beta}: the noises are Gaussian white and independent with intensities~$\sigma_i$, $i=v,\omega,a,\beta$, and the biases are constant; the numerical values are given in table~\ref{tab:meas}. The tuning matrices are set to $(K,L,M):=(5I,5I,0.5I)$.
\begin{table}[t!]
	\setlength{\extrarowheight}{2pt}
	\centering
	\begin{tabular}{| l | l || l | l |}
		\firsthline
		$b_v$ & $(-0.10~0.30~-0.05)^T$ & $\sigma_v^2$ & $\num{e-5}(2~2~2)^T$ \\ \hline
		$b_\omega$ & $(0.0250~-0.0300~-0.0175)^T$ & $\sigma_\omega^2$ & $\num{e-7}(2~2~2)^T$ \\ \hline
		$b_a$ &  $(0.05~0.04~-0.02)^T$ & $\sigma_a^2$ & $\num{e-5}(1~1~1)^T$ \\ \hline
		$b_\beta$ &  $(0.024 -0.020 -0.018)^T$ & $\sigma_\beta^2$ & $\num{e-7}(1~1~1)^T$ \\
		\lasthline
	\end{tabular}
	\caption{Biases and noises in the actual measurements.\label{tab:meas}}
\end{table}

\begin{figure}[ht]\hspace{-5mm}
	\includegraphics[width=1.1\columnwidth]{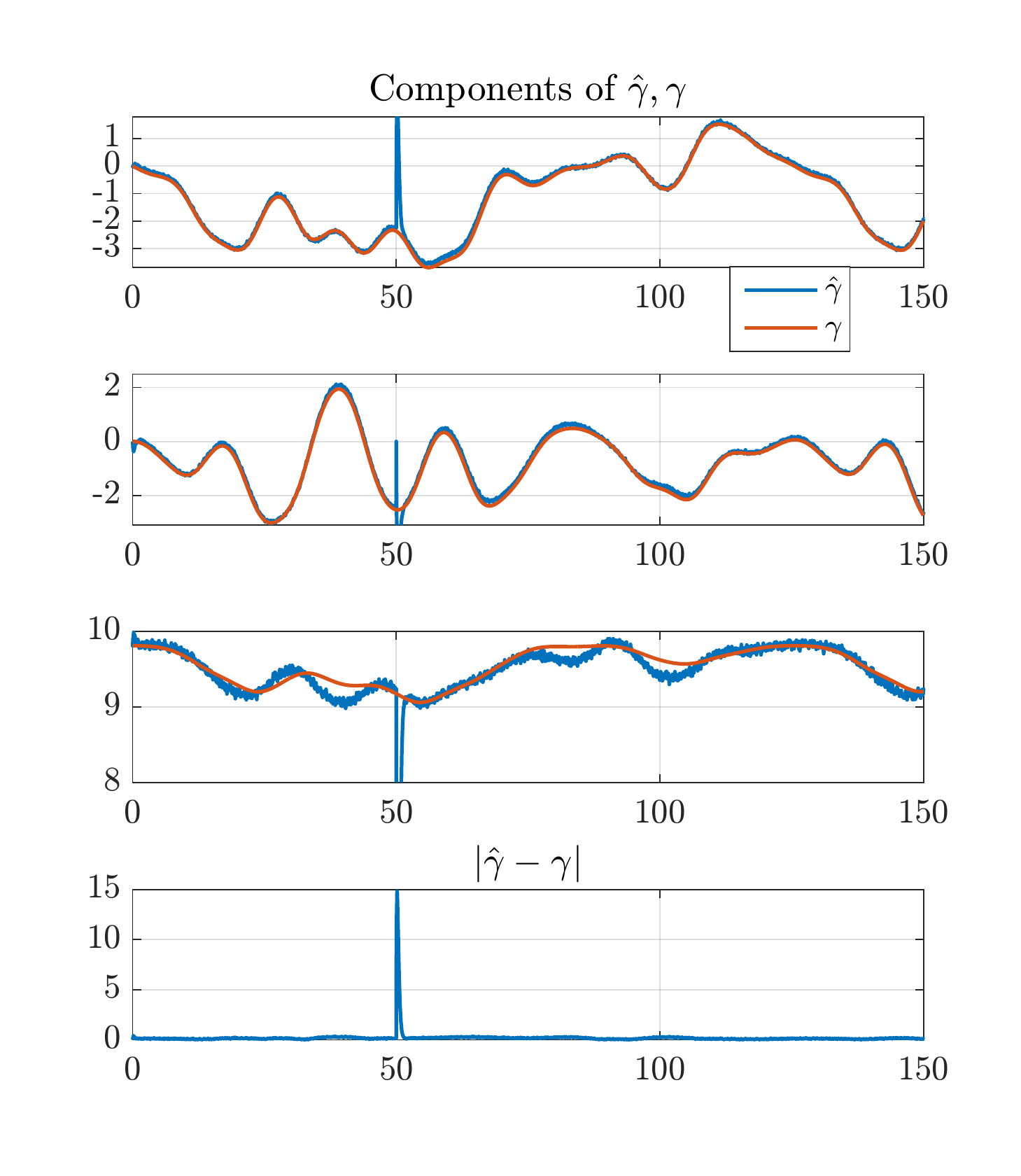}
	\caption{Components of true $\gamma$ (red), measured $\gamma_m$ (blue) and estimated~$\hat\gamma$ (orange).}
	\label{fig:gamma}
\end{figure}

The observer is initialized without error, but is suddenly reinitialized at~$t=50$. The convergence of the estimated states is as anticipated excellent after the reinitialization, very fast for $\hat v,\hat\gamma$ and slower for $\hat\beta$, in accordance with the choice of gains, see Fig.~\ref{fig:v}, \ref{fig:gamma} and~\ref{fig:beta}. The (desirable) independence between $\hat v,\hat\gamma$ and~$\hat\beta$ is clearly visible for $t\in[80,100]$: only $\hat\beta$ is affected by the disturbance of the magnetic field~$B$.

\begin{figure}[ht]\hspace{-5mm}
	\includegraphics[width=1.1\columnwidth]{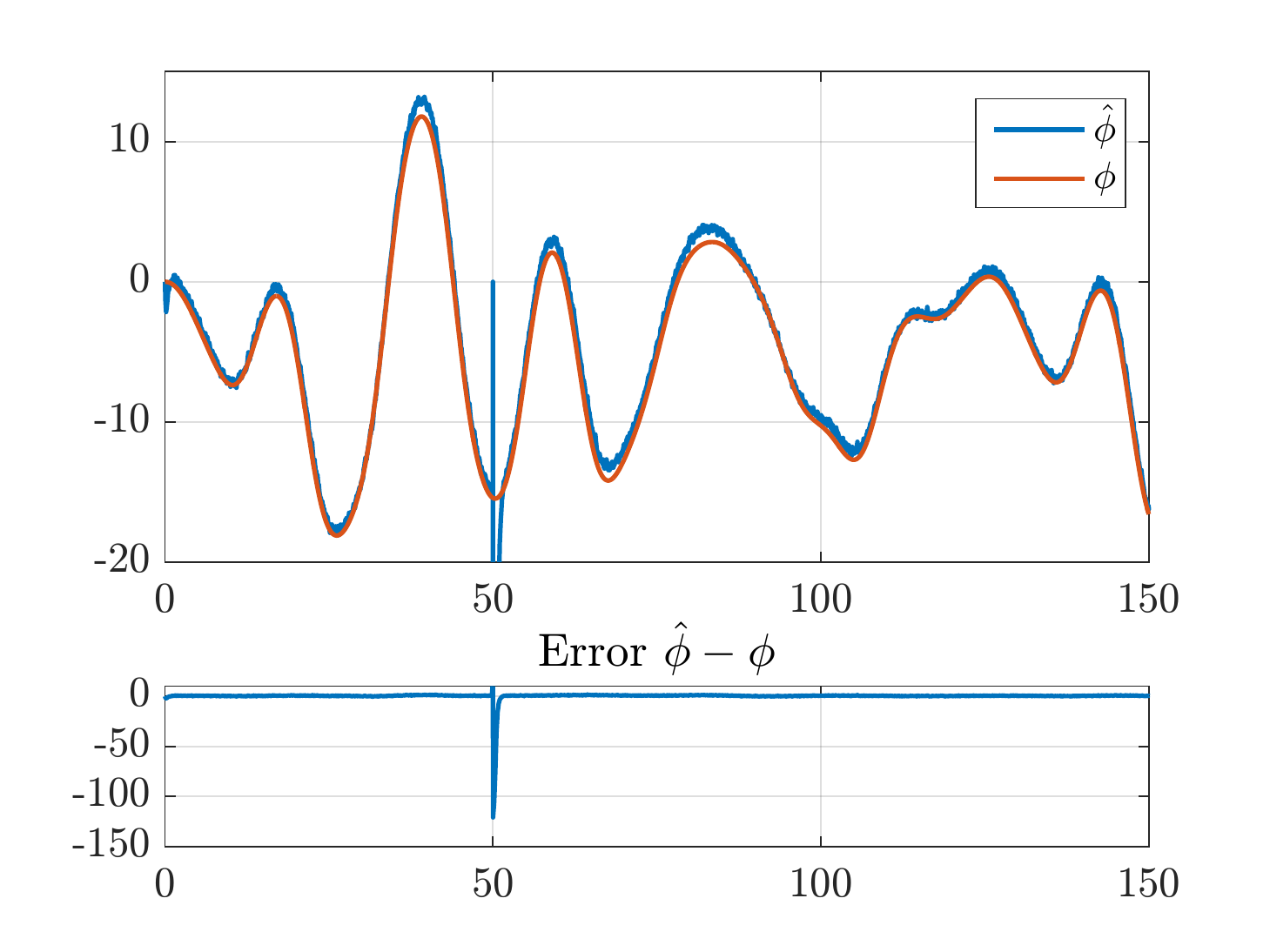}
	\caption{True $\phi$ (red) and estimated $\hat\phi$ (blue).}
	\label{fig:phi}
\end{figure}

\begin{figure}[ht]\hspace{-5mm}
	\includegraphics[width=1.1\columnwidth]{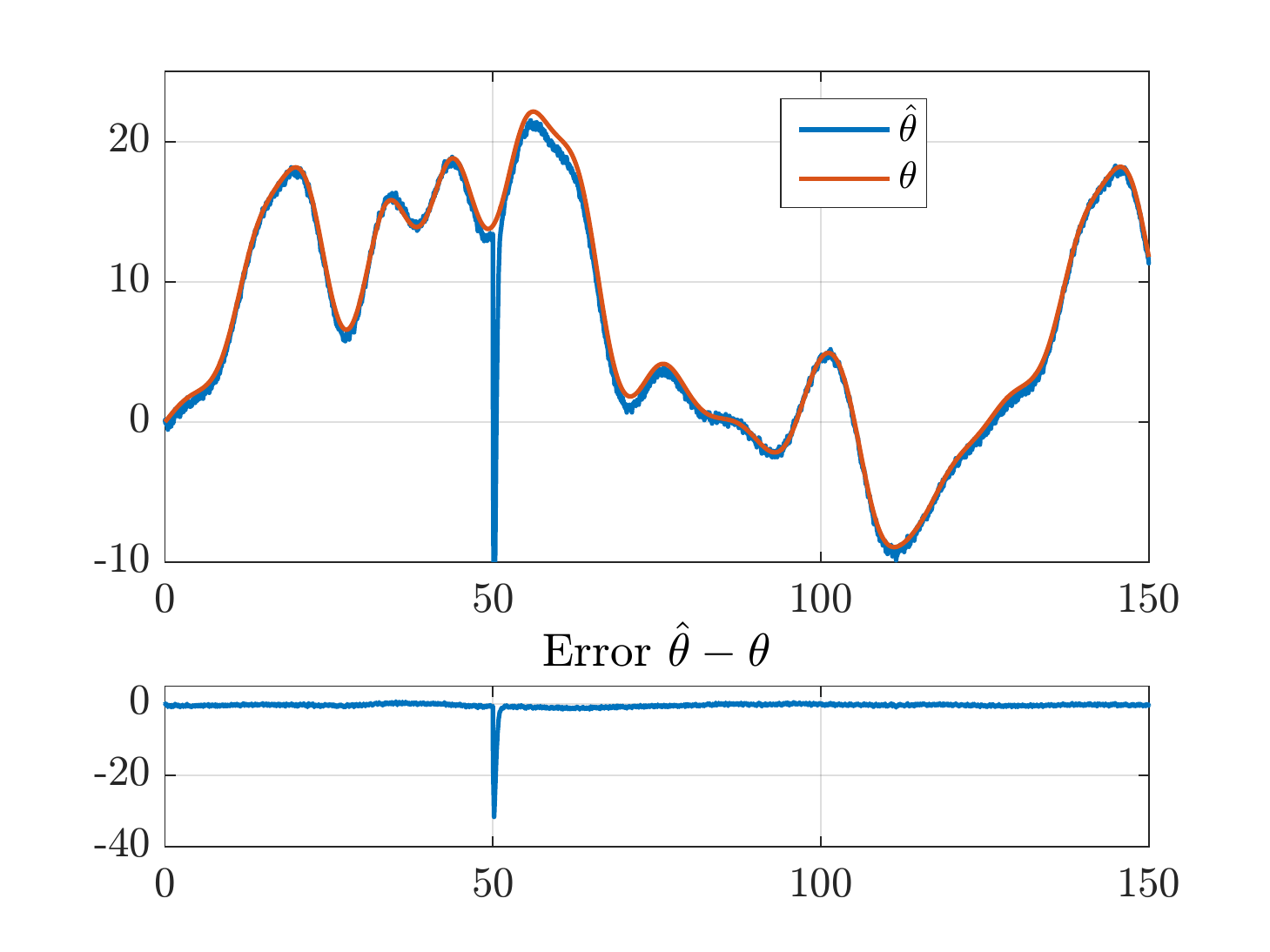}
	\caption{True $\theta$ (red) and estimated $\hat\theta$ (blue).}
	\label{fig:theta}
\end{figure}

\begin{figure}[ht]\hspace{-5mm}
	\includegraphics[width=1.1\columnwidth]{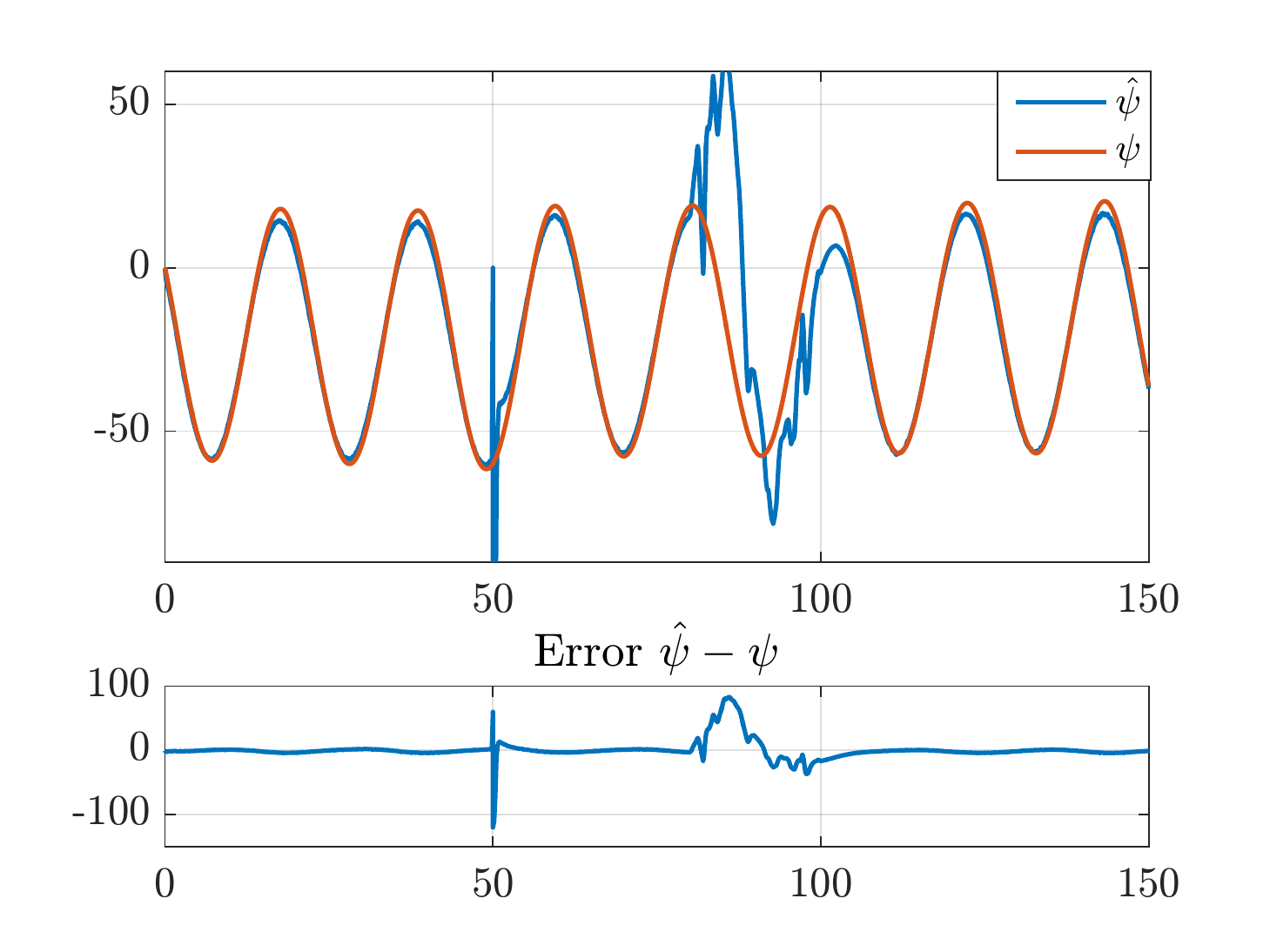}
	\caption{True $\psi$ (red) and estimated $\hat\psi$ (blue).}
	\label{fig:psi}
\end{figure}

Finally, Fig.~\ref{fig:phi}--\ref{fig:psi} show the Euler angles $\phi,\theta,\psi$ (in degrees) corresponding to the estimated orientation reconstructed from $\hat\gamma,\hat\beta$, using the nominal values of $g$ and~$B$. Notice the pitch angle~$\theta$ and roll angle~$\phi$ are as anticipated unaffected by the magnetic disturbance for $t\in[80,100]$.

We emphasize that the small discrepancies between the true states and their estimates once the observer has converged are due only to the imperfect measurements and to the magnetic disturbance. Without this effects, which are not taken into account in the model, the estimated and true states perfectly agree, as predicted by the theoretical analysis.

\section{Conclusion}
We have presented a simple ``geometry-free'' observer for estimating the attitude and velocity of a rigid body from the measurements of specific acceleration, angular velocity, magnetic field (in body axes), and linear velocity (in body axes). We have established its global exponential convergence by a very simple yet rigorous Lyapunov analysis. This is an improvement the existing literature, where only the almost global asymptotic convergence is achieved, moreover at the cost of a more involved analysis.

The simple structure of the observer and its strong convergence properties are promising for tackling the case where measurement biases are explicitly considered, as well as for designing closed-loop controllers with output feedback.


\bibliographystyle{plain}        
\bibliography{ref}           





\end{document}